\documentclass[10pt,openbib]{article}

\newcommand{\tsum}{\tilde\#_S}
\newcommand{\mc}[1]{\mathcal{#1}}
\newcommand{\ii}{^{(i)}}

\newcommand{\limfunc}[1]{\textnormal {#1\,}}
\newcommand{\func}[1]{\textnormal {#1\,}}

\newcommand{\tublim}{\overset{C^\infty}{\underset{S\rightarrow\infty}{\longrightarrow}}}

\newcommand{\hol}{\func{Hol}}
\newcommand{\img}{\func{img}}
\newcommand{\tr}{\func{tr}}
\newcommand{\e}[1]{\mathbf{e}^{#1}}
\newcommand{\R}{\mathbb{R}}

\newcommand{\vol}{\func{Vol}}

\usepackage[UKenglish]{babel}
\usepackage{amssymb}
\usepackage{makeidx}
\usepackage{amsfonts}
\usepackage{amsmath}
\usepackage{graphicx}
\usepackage{hyperref}
\usepackage[all]{xy}

\newtheorem{theorem}{Theorem}

\newtheorem{claim}[theorem]{Claim}

\newtheorem{condition}[theorem]{Condition}
\newtheorem{conjecture}[theorem]{Conjecture}
\newtheorem{corollary}[theorem]{Corollary}

\newtheorem{definition}[theorem]{Definition}
\newtheorem{example}[theorem]{Example}

\newtheorem{lemma}[theorem]{Lemma}
\newtheorem{notation}[theorem]{Notation}

\newtheorem{proposition}[theorem]{Proposition}

\newenvironment{proof}[1][Proof]{\noindent 
         \begin{list}{}{\setlength{\topsep}{0cm} \setlength{\leftmargin}{0.6cm}}
          \item[\hspace{0.19cm}\textbf{#1} \hspace{0.4cm}]}{  \rule{0.5em}{0.5em} \end{list}}

\begin{document}

%+Title
\title{$G_2-$instantons over Kovalev manifolds II}
\author{Henrique N. S\'a Earp}
\date{\today}
\maketitle
%-Title

%+Abstract
\begin{abstract}
This is the first nontrivial construction to date of instantons over a compact manifold with $\func{Hol}=G_2$. The HYM connections on asymptotically stable bundles over Kovalev's noncompact Calabi-Yau 3-folds, obtained in the first article \cite{G2I}, are glued compatibly with a twisted connected sum, to produce a $G_2-$instanton over the resulting compact $7-$manifold \cite{kovalevzinho,kovalevzao}.
This is accomplished under a nondegeneracy \emph{acyclic} assumption on the bundle `at infinity', which occurs e.g. over certain projective varieties $X_{22}\hookrightarrow \mathbb{C}P^{13}$ \cite{X22,Mukai Fano 3-folds} equipped with an  \emph{asymptotically rigid} bundle.
\end{abstract}
%-Abstract

%+Contents
\tableofcontents
%-Contents

\setcounter{section}{0}
\section*{Introduction}
\addcontentsline{toc}{section}{Introduction and acknowledgements}

The first paper \cite{G2I} was devoted to the Hermitian Yang-Mills (HYM) problem over A. Kovalev's noncompact Calabi-Yau 3-folds $W$  with an exponentially asymptotically cylindrical (EAC) end.
Its main theorem asserts the existence of such a metric on holomorphic bundles $\mathcal{E}\rightarrow W$ which are \emph{asymptotically stable}, i.e., (slope-)stable over the divisor `at infinity' $D\subset \bar{W}$. Moreover, solutions have  `$C^\infty-$exponential decay'  along the tubular end to a reference metric, which extends the instanton metric on $\left.\mc{E}\right\vert_D$ in a prescribed way. On the other hand, it is known that a HYM connection $A$ on a holomorphic vector bundle    over a Calabi-Yau $3-$fold $\left(W,\omega,\Omega\right)$ lifts to a $G_{2}-$instanton
on the pull-back bundle $p_1^*\mathcal{E}\rightarrow M=W\times S^{1}$, hence the above result yields a non-trivial solution of the corresponding instanton equation:
\begin{equation}        \label{sd part kahler -> G2}
        F_A\wedge\ast\varphi
        \doteq\tfrac{1}{2}F_{A}\wedge \left( \omega \wedge \omega -2\func{Re}%
        \Omega \wedge d\theta \right)=0  
\end{equation}
(where $d\theta$ is the coordinate 1-form on $S^1$ and $\varphi$ is the induced $G_2-$structure).

This paper will be concerned with the natural sequel, the gluing of two such solutions according to Kovalev's twisted connected sum of the base manifolds. That construction joins two `truncated' EAC products $W_{S}^{(i)}\times S^{1}$ in order to obtain a smooth \emph{compact} $G_2-$manifold
\begin{equation*}
M_{S}=\left( W_{S}^{\prime }\times S^{1}\right) \cup _{F_{S}}\left(
W_{S}^{\prime \prime }\times S^{1}\right) \doteq W'\widetilde{\#}_{S}W'',
\end{equation*}%
where $S$ is the `neck-length' parameter and $F_S$ is the product of a hyper-K\"ahler rotation on the divisors `near infinity' and a nontrivial identification of the circle components of the boundary. In view of the exponential decay property of solutions over each end, one is led to expect that residual self-dual curvature, i.e., error terms from truncation in the corresponding $G_2-$instanton equation over $M_S$, can be dealt with by a perturbative, `stretch-the-neck'-type argument. 

This paper's main result [\emph{Theorem \ref{thm gluing}}, \emph{Subsection \ref{subsec Statement of gluing theorem}}] posits that this is indeed the case, albeit under a rigidity assumption on the `bundle at infinity'. In order to guarantee a right-inverse for the linearisation of the self-dual curvature operator $A\mapsto F_A\wedge\ast\varphi$ around a solution, one requires that the corresponding deformation complex over $\left.\mc{E}\right\vert_D$ be \emph{acyclic}. In other words, the instanton `at infinity' must be an isolated point in its moduli space. Examples satisfying this hypothesis are obtained from certain base manifolds $X_{22}\hookrightarrow\mathbb{C}P^{13}$ admissible by Kovalev's construction and studied by Iskovskih \cite{X22} and Mukai \cite{Mukai Fano 3-folds}. 

Readers familiar with the 4-dimensional model outlined by Donaldson in \cite{floer}
 will find that this paper mimics that source in all its essential aspects.

\setcounter{section}{0}
\section*{Acknowledgements}

This second paper was written at Unicamp, supported by post-doctoral grant 2009/10067-0 from Fapesp - S\~ao Paulo State Research Council. It contains some unpublished material from my thesis, which was funded by Capes Scholarship 2327-03-1, from the Brazilian Ministry of Education.

I would like to thank Simon Donaldson, for suggesting this outcome of the first article and guiding its initial steps, and Marcos Jardim, for many useful subsequent discussions.

\section{$G_2-$instantons over EAC manifolds}

Let us briefly recall the basic facts about the Calabi-Yau $3-$folds $%
W_{i} $ that one intends to glue together, via a nontrivial product with $S^1$.
\label{def base manifold (W,w)}
\index{base manifold}%
A \emph{base manifold} for our purposes is a compact, simply-connected K\"{a}hler $3-$fold $\left( \bar{W},\bar{\omega}\right) $ carrying a $K3-$divisor  $D \in \left\vert -K_{\bar{W}}\right\vert$ with holomorphically trivial         normal bundle $\mathcal{N}_{D/\bar{W}}$ such that the complement  $W\doteq\bar{W}\setminus D$ has $\pi _{1}\left( W\right) $ finite.
Topologically $W$ looks like a compact manifold $W_{0}$ with boundary $%
D\times S^{1}$ and a cylindrical end attached:%
\begin{equation}        \label{eq cylindrical picture}
\begin{array}{c}
        W=W_{0}\cup W_{\infty } \\ 
        W_{\infty }\simeq  D\times S_\alpha^{1}
        \times \left(\mathbb{R}_{+}\right)_s .%
\end{array}
\end{equation}

The $K3$ divisor $D$ is actually hyper-K\"ahler, with complex
structure $I$ inherited from $\bar{W}$ and the additional structures $J$ and $K=IJ$ satisfying the quaternionic
relations; denote  by $\kappa_I,%
\kappa _{J}$ and $\kappa _{K}$ their associated K\"{a}hler forms. Then \cite[Theorem 2.2]{kovalevzinho} $W$ admits a complete Calabi-Yau structure $\omega $ and holomorphic volume form $\Omega $ \emph{exponentially asymptotic}  to the cylindrical model
\begin{eqnarray}        
\index{Kahler metric@K\"{a}hler metric!winf@$\omega _{\infty }$}
        \begin{array}{r c l}      \label{asymptotic forms}
                \omega _{\infty } &=&\kappa _{I}+ds\wedge d\alpha \\
                \Omega _{\infty } &=&\left( ds+\mathbf{i}d\alpha \right)                 \wedge \left(\kappa _{J}+\mathbf{i}\kappa _{K}\right) \text{,}
        \end{array}
        \end{eqnarray}
in the sense that
\begin{displaymath}
        \left.\omega \right\vert_{W_{\infty}}
        =\omega _{\infty }+d\psi, \qquad
        \left.\Omega \right \vert_{W_\infty}
        =\Omega _{\infty }+d\Psi, 
\end{displaymath}%
where  $\psi $ and  $\Psi $ are smooth 
and decay exponentially in all derivatives along the tubular end.

As to the gauge-theoretic initial data \cite{G2I}, let $z=\e{-s+\bf{i}\alpha} $ be the holomorphic coordinate along the tube and denote $D_z$ the corresponding $K3$ component of the boundary. 
\label{def bundle E->W}
\index{infinity!stable at} 
A bundle $\mathcal{E}\rightarrow W$ is called \emph{asymptotically stable} (or \emph{stable at infinity}) if it is the
restriction of an indecomposable holomorphic vector bundle $\mathcal{E}\rightarrow \bar{W}$ such that
        $\left. \mathcal{E}\right\vert _{D}$ is stable (hence also $\left. \mathcal{E}\right\vert _{D_{z}}$ for $\left\vert z\right\vert <\delta $).
Moreover,
\label{def reference metric H0}
\index{metric!reference $H_0$}
a \emph{reference metric} $H_{0}$ on such 
$\mathcal{E}\rightarrow W$ is (the restriction of) a smooth Hermitian metric on $\mathcal{E}\rightarrow \bar{W}$ such that $\left. H_{0}\right\vert _{D_{z}}$ are the corresponding HYM metrics on $\left. \mathcal{E}\right\vert _{D_z}$, $0\leq \left\vert z\right\vert <\delta$, and $H_{0}$ has finite energy: 
$\Vert\hat{F}_{H_{0}}\Vert _{L^{2}\left( W,\omega\right) }<\infty$.

Then, given an asymptotically stable bundle with reference metric $\left(\mc{E},H_0  \right)$, a nontrivial smooth $G_2-$instanton on $p_1^*\mc{E}\rightarrow W\times S^1$ is obtained from every solution of the HYM problem over $W$. Moreover, such solutions have the property of exponential asymptotic decay in all derivatives to $H_0$ along the tubular end [\emph{ibid.}, Theorem 59]: %
\index{asymptotic decay!of HYM metric}%
\begin{equation}       \index{metric!Hermitian Yang-Mills}
                        \label{eq HYM solution}
\fbox{$\begin{array}{c}
        \hat{F}_{H}=0, \quad H \tublim H_0. \\
\end{array}$}
\end{equation}
Here convergence takes place over cylindrical bands of fixed `size':

\begin{notation}        \label{Not finite cylinder}%
Let $Q\rightarrow W$ be a bundle equipped with a fibrewise metric and denote $W_S$ the truncation of $W$ at `neck length' $S$; given $S>r>0$, write $\Sigma _{r}\left( S\right) $ for the
interior of the cylinder $\left( W_{S+r}\smallsetminus
W_{S-r}\right) $ of `length' $2r$. We denote the \emph{$C^k-$exponential tubular limit} of an element in  $C^k\left(\Gamma(Q)\right)$ by:
\begin{eqnarray*}
        \phi \overset{C^k}{\underset{S\rightarrow\infty}{\longrightarrow}}\phi_{0}         & \dot{\Leftrightarrow}& \left\Vert \phi -\phi_{0} \right\Vert_{C^k\left( \Sigma_1(S),\omega \right)}=O\left( \e{-S} \right).
\end{eqnarray*}
\end{notation}

Finally, 
let us fix some vocabulary towards the statement of the main theorem. Denote $A_{0}$ the Chern connection of $H_0$; then by definition each $\left. A_{0}\right\vert _{D_{z}}
$ is ASD. In particular, $\left.
A_{0}\right\vert _{D}$ induces an elliptic deformation complex%
\begin{equation}        \label{eq ell complex over D}
\Omega ^{0}\left(\left.\mathfrak{g}\right\vert_D\right) \overset{d_{A_{0}}}{\rightarrow }%
\Omega ^{1}\left( \left.\mathfrak{g}\right\vert_D\right) \overset{d_{A_{0}}^{+}}{\rightarrow }%
\Omega _{+}^{2}\left( \left.\mathfrak{g}\right\vert_D\right) 
\end{equation}%
where $\left.\mathfrak{g}\right\vert_D=%
\limfunc{Lie}\left( \left. \mathcal{G}\right\vert _{D}\right) $ generates
the gauge group $\mathcal{G}$ $=\limfunc{End}\mathcal{E}$ over $D$. Thus,
the requirement that $\mathcal{E}$ be indecomposable restricts the
associated cohomology: 
\begin{equation*}
        \mathbf{H}_{A_{0}\vert_D}^{0}=0.
\end{equation*}%
On the other hand, one might restrict
attention to \emph{acyclic }connections,
i.e., whose gauge class $\left[ A_{0}\right] $ is \emph{isolated} in $%
\mc{M}_D\doteq\mc{M}_{\left. \mathcal{E}\right\vert _{D}}$. The absence of  infinitesimal deformations translates into the vanishing of the other cohomology
group:%
\begin{equation*}
        \mathbf{H}_{A_{0}\vert_D}^{1}=0.
\end{equation*}

\begin{definition}      \label{def asymp rigid}
A reference metric on  $\mathcal{E}\rightarrow W$ is \emph{asymptotically rigid} if the associated complex $(\ref{eq ell complex over D})$ over $D$ has trivial cohomology.
\end{definition}

\subsection{Suitable pairs and gluing}

A $7-$dimensional product $W_{S}^{\prime }\times S^{1}$, where $W'$ is of the above form, has boundary $%
D^{\prime }\times S^{1}\times S^{1}$. Comparing the asymptotic model $\left( %
\ref{asymptotic forms}\right) $ with the standard form of the $G_{2}-$%
structure on a product $CY\times S^{1}$ we find that $W_{S}^{\prime }\times S^{1}$ carries a $G_{2}-$structure on a collar neighbourhood of the boundary that is asymptotic to:%
\begin{equation}        \index{G2@$G_{2}$!-structure@$-$structure}
                        \label{asymptotic G2-structure}
        \varphi _{S}^{\prime }=\kappa _{I}^{\prime }\wedge d\alpha 
        +\kappa_{J}^{\prime }\wedge d\theta +\kappa _{K}^{\prime }
        \wedge ds+d\alpha \wedge d\theta \wedge ds.
\end{equation}%
Since the inclusion of the set of all $G_{2}-$structures $\mathcal{P}^{3}\left( W^{\prime
}\times S^{1}\right) \subset \Omega ^{3}\left( W^{\prime }\times
S^{1}\right) $ is open \cite[p. 243]{Joyce}, $\varphi _{T}^{\prime }$ is
itself a $G_{2}-$structure on $W^{\prime }\times S^{1}$ for large $S$.
\begin{condition}
Two manifolds $W^{\prime }$ and $W^{\prime \prime }$ as above will be
suitable for the gluing procedure if there is a hyper-K\"{a}hler isometry \begin{equation*}
        f:D_{J}^{\prime }\rightarrow D^{\prime \prime }
\end{equation*}%
between $D^{\prime \prime }$ and the hyper-K\"{a}hler rotation of $D^{\prime
}$ with complex structure $J$. In this case the
(pull-back) action on K\"{a}hler forms is%
\begin{equation}        \index{K3 surface!hyper-K\"{a}hler}%
f^{\ast }:\kappa _{I}^{\prime \prime }\mapsto \kappa _{J}^{\prime },\qquad
\kappa _{J}^{\prime \prime }\mapsto \kappa _{I}^{\prime },\qquad \kappa
_{K}^{\prime \prime }\mapsto -\kappa _{K}^{\prime }.
\label{matching of kahler forms}
\end{equation}
\end{condition}

Assuming this holds, define a map between collar neighbourhoods of the
boundaries by 
\begin{eqnarray*}
        F_{S}:D^{\prime }\times S^{1}\times S^{1}\times \left[ S-1,S\right]
        &\rightarrow& 
        D^{\prime \prime }\times S^{1}\times S^{1}\times \left[ S-1,S\right]         \\
        \left( y,\alpha ,\theta ,s\right) 
        &\mapsto& 
        \left( f\left( y\right) ,\theta,\alpha ,2S-1-s\right) .
\end{eqnarray*}
This identification gives a compact oriented $7-$manifold 
\begin{equation*}
M_{S}=\left( W_{S}^{\prime }\times S^{1}\right) \cup _{F_{S}}\left(
W_{S}^{\prime \prime }\times S^{1}\right) \doteq W'\widetilde{\#}_{S}W''.
\end{equation*}%

\begin{figure}[h]
        
%MINHAS FIGURAS

\setlength{\unitlength}{0.5cm}

\begin{picture}(0,8.5)(1,-3.5)

%%%%% INICIO FIGURA ESQUERDA

%ORIGEM
%\put(0,0){\circle*{0.2}}\put(-1,-0.5){(0,0)}

%QUADRICULADO
%\multiput(0,-5)(1,0){25}{\line(0,1){12}}
%\multiput(0,-5)(0,1){13}{\line(1,0){18}}

\put(1.3,-0.5){\huge{$W'_S$}}

\linethickness{0.4mm}

%CURVA CONTINUA
%%CIMA
\qbezier(4,2)   (3,3)
        (2,3)

\qbezier(4,2)   (5,1)
        (6,1)

%%BAIXO
\qbezier(2,-3)  (3,-3)
        (4,-2)
\qbezier(4,-2)   (5,-1)
        (6,-1)

%LATERAL
\qbezier(2,3)(0,3)
        (0,0)
\qbezier(0,0)(0,-3)
        (2,-3)

%TUBO
\put(9,1){\line(-1,0){3}}
\put(9,-1){\line(-1,0){3}}

%SECAO TRANSVERSAL

\qbezier(8,1)  (7,0)
        (8,-1)

\linethickness{0.7mm}

\qbezier(9,1)  (8,0)
        (9,-1)

\qbezier(9,1)  (10,0)
        (9,-1)

%VETORES
\thicklines
\put(7.5,0){\line(1,0){1}}
\put(7.6,0){\vector(1,0){0,7}}

\put(7.6,0.5){\line(1,0){1}}
\put(7.7,0.5){\vector(1,0){0,7}}

\put(7.6,-0.5){\line(1,0){1}}
\put(7.7,-0.5){\vector(1,0){0,7}}

%SOMA CONEXA TORCIDA

\put(9.2,0.8){
\begin{xy}(0,0);<0.34cm,0cm>:
%%SUPERFICIE D'
        (-0.3,1.5)*{\txt\small\itshape{D}}="D'",*!/^2.5mm/{'_J},
        (-1,0);(1,0) **\crv{(0,1)},
        (-1,0);(0,-0.5) **\crv{(0,0)},
        (0,-0.5);(1,0) **\crv{(0,0)},
        (0.04,-1.15) *{^\times},
        %
%%S1 ALPHA'
        (-1.75,-1.7) *{\txt\small\itshape{S}},p+(-0.9,-1.25)*{^1_{\alpha'}},
        (0,-1.5);(0,-2.5) **\crv{(0.5,-2)},
        ?(0.9)*!/^0.2mm/{\txt{ \^}},
        (0,-1.5);(0,-2.5) **\crv{(-0.5,-2)},
        (0.04,-3.2) *{^\times},
        %
%%S1 THETA'
        (-1.75,-4.3) *{\txt\small\itshape{S}},p+(-0.9,-2.75)*{^1_{\theta'}},
        (0,-3.5);(0,-4.5) **\crv{(0.5,-4)},
        ?(0.9)*!/^0.2mm/{\txt{ \^}},
        (0,-3.5);(0,-4.5) **\crv{(-0.5,-4)},
        %
%%ZOOM NO BORDO    
        (-3,-2.8);p+(2,0) **@{~},
        (-3,0.4);p+(2,0) **@{~},
        %
%IDENTIFICACOES
        (0,-1.5);(4,-3.5) **@{.},
        (0,-2,5);(4,-4.5) **@{.},
        (0,-3.5);(4,-1.5) **@{.},
        (0,-4,5);(4,-2.5) **@{.},
        (1.2,0);(2.8,0) **@{.} 
        ?<*@{<} ?(0.5)*!/_3mm/\txt\small\itshape{f}?>*@{>},
        "D'"+(-4.5,1.5);p+(14.5,0) **\frm{^)}
        ?(0.5)*!/_15mm/{\cup_{F_S}}, 
\end{xy}
}

\put(13.3,0.8){
\begin{xy}(0,0);<0.34cm,0cm>:
%%SUPERFICIE D''
        (0,1.5)*{\txt\small\itshape{D}},*!/^2.5mm/{''},
        (-1,0);(1,0) **\crv{(0,1)},
        (-1,0);(0,-0.5) **\crv{(0,0)},
        (0,-0.5);(1,0) **\crv{(0,0)},
        (0.04,-1.15) *{^\times},
        %
%%S1 ALPHA''
        (1,-1.7) *{\txt\small\itshape{S}},p+(2,-1.25)*{^1_{\alpha''}},
        (0,-1.5);(0,-2.5) **\crv{(0.5,-2)},
        ?(0.9)*!/^0.2mm/{\txt{ \^}},
        (0,-1.5);(0,-2.5) **\crv{(-0.5,-2)},
        (0.04,-3.2) *{^\times},
        %
%%S1 THETA''
        (1,-4.3) *{\txt\small\itshape{S}},p+(2,-2.75)*{^1_{\theta''}},
        (0,-3.5);(0,-4.5) **\crv{(0.5,-4)},
        ?(0.9)*!/^0.2mm/{\txt{ \^}},
        (0,-3.5);(0,-4.5) **\crv{(-0.5,-4)},
        %
%%ZOOM NO BORDO    
        (0.9,-2.8);p+(2,0) **@{~},
        (0.9,0.4);p+(2,0) **@{~},
\end{xy}
}

%%%%% FIGURA DIREITA %%%%%
%
\begin{picture}(10,10)(-15.7,0)

%ORIGEM
%\put(0,0){\circle*{0.2}}\put(-1,-0.5){(0,0)}

%QUADRICULADO
%\multiput(0,-5)(1,0){13}{\line(0,1){12}}
%\multiput(-5,-5)(0,1){13}{\line(1,0){10}}

\put(6.2,-0.5){\huge{$W''_S$}}

\linethickness{0.4mm}

%CURVA CONTINUA
%%CIMA
\qbezier(6,2)   (7,3)
        (8,3)

\qbezier(6,2)   (5,1)
        (4,1)

%%BAIXO
\qbezier(8,-3)  (7,-3)
        (6,-2)
\qbezier(6,-2)   (5,-1)
        (4,-1)

%LATERAL
\qbezier(8,3)(10,3)
        (10,0)
\qbezier(10,0)(10,-3)
        (8,-3)

%TUBO
\put(1,1){\line(1,0){3}}
\put(1,-1){\line(1,0){3}}

%SECAO TRANSVERSAL
\qbezier(2,1)  (1,0)
        (2,-1)

\thinlines

\qbezier(2,1)  (3,0)
        (2,-1)

\linethickness{0.7mm}

\qbezier(1,1)  (0,0)
        (1,-1)

\thicklines
\put(0.5,0){\line(1,0){1}}
\put(0.6,0){\vector(1,0){0,7}}

\put(0.6,0.5){\line(1,0){1}}
\put(0.7,0.5){\vector(1,0){0,7}}

\put(0.6,-0.5){\line(1,0){1}}
\put(0.7,-0.5){\vector(1,0){0,7}}

\end{picture}
%%%%% FIM DA FIGURA DIREITA %%%%

\end{picture}
        \caption{The compact $7-$manifold $M_S$}
\end{figure}
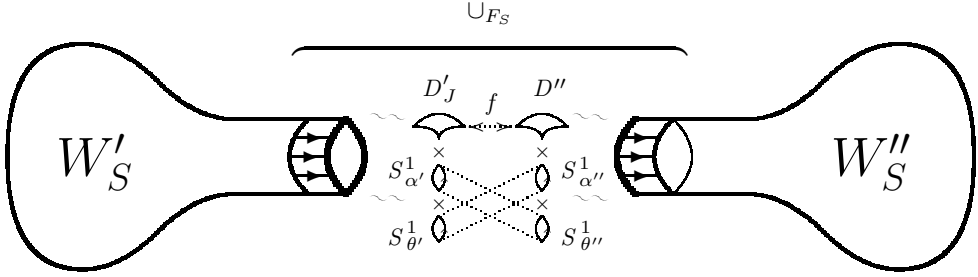

The matching of K\"{a}hler forms $\left( \ref{matching of kahler forms}%
\right) $ guarantees that the respective $G_{2}-$structures $\left( \ref%
{asymptotic G2-structure}\right) $ on $W_{S}^{\prime }\times S^{1}$ and $%
W_{S}^{\prime \prime }\times S^{1}$ agree along the gluing region $\left[
S-1,S\right] $:%
\begin{eqnarray*}
F_{S}^{\ast }\varphi _{S}^{\prime \prime } &=&F_{S}^{\ast }(\underset{\kappa
_{J}^{\prime }}{%
\underbrace{\kappa _{I}^{\prime \prime }}}\wedge \underset{d\theta ^{\prime }%
}{\underbrace{d\alpha ^{\prime \prime }}}+\underset{\kappa _{I}^{\prime }}{%
\underbrace{\kappa _{J}^{\prime \prime }}}\wedge \underset{d\alpha ^{\prime }%
}{\underbrace{d\theta ^{\prime \prime }}}+\underset{-\kappa _{K}^{\prime }}{%
\underbrace{\kappa _{K}^{\prime \prime }}}\wedge \underset{-ds^{\prime }}{%
\underbrace{ds^{\prime \prime }}}+\underset{d\theta ^{\prime }}{\underbrace{%
d\alpha ^{\prime \prime }}}\wedge \underset{d\alpha ^{\prime }}{\underbrace{%
d\theta ^{\prime \prime }}}\wedge \underset{-ds^{\prime }}{\underbrace{%
ds^{\prime \prime }}}) \\
&=&\kappa _{I}^{\prime }\wedge d\alpha ^{\prime }+\kappa _{J}^{\prime
}\wedge d\theta ^{\prime }+\kappa _{K}^{\prime }\wedge ds^{\prime }+d\alpha
^{\prime }\wedge d\theta ^{\prime }\wedge ds^{\prime } \\
&=&\varphi _{S}^{\prime }.
\end{eqnarray*}%
so we obtain a globally well-defined $G_{2}-$structure $\varphi _{S}$ on $%
M_{S}$. Thus, for large enough $S$, there is a $1-$parameter family $\left(
M_{S},\varphi _{S}\right) $ of compact oriented manifolds $M_{S}$ equipped
with $G_{2}-$structures $\varphi _{S}$. 

While it is possible
to arrange $d\varphi _{S}=0$ for any $S$ \cite[eq. $\left( 4.23\right) $]%
{kovalevzao}, a pair $\left( M_{S},\varphi _{S}\right) $ is not in principle
a $G_{2}-$manifold, as one has yet to satisfy the co-closedness
condition:%
\begin{equation*}
        d\ast _{\varphi _{S}}\varphi _{S}=0.
\end{equation*}
In fact, although the cut-off functions involved in the asymptotic approximations leading to $\left( \ref{asymptotic G2-structure}\right) $   add error terms to $%
d\ast _{\varphi _{S}}\varphi _{S}$, these are controlled by the estimate \cite[\emph{Lemma 4.25}]{kovalevzao}
\begin{equation*}
        \left\Vert d\ast _{\varphi _{S}}\varphi _{S}\right\Vert _{L_{k}^{p}}
        \leq C_{p,k}e^{-\lambda S},
\end{equation*}%
with $0<\lambda <1$. This exponential decay implies that, by `stretching
the neck' up to large enough $S_{0}$, one can make the error so small as to be compensated by a suitably small perturbation of $\varphi _{S}$ in 
$\mathcal{P}^{3}\left( M_{S}\right) $, $S>S_{0}$ \cite[Theorem 2.3]%
{kovalevzinho}. Hence one achieves a $1-$parameter family of compact oriented $%
G_{2}-$manifolds:%
\begin{equation*}       \index{G2@$G_{2}$!-manifold@$-$manifold}
        \left( M_{S},\tilde \varphi _{S}\right) ,\qquad S>S_{0}.
\end{equation*}

\subsection{Statement of the instanton gluing theorem}
\label{subsec Statement of gluing theorem}

Let $\left(M_S,\tilde\varphi_S\right)$ be a compact  $G_2-$manifold   with $\hol(\tilde\varphi_S)= G_2$ as above, obtained from a Fano pair by Kovalev's construction:
\begin{displaymath}
        M_S = W'\tsum W'' \doteq \left( W'_S \times S^1 \right)
        \cup_{F_S} \left( W''_S \times S^1 \right). 
\end{displaymath}
Furthermore, let $\mc{E}^{(i)}\rightarrow W^{(i)}$ be asymptotically stable holomorphic bundles with same structure group $G=\func{Aut}(\mc{E})$, such that there is a $G-$isomorphism \begin{displaymath}
        g:\left.\mc{E}'\right\vert_{D_J'}\tilde\rightarrow \left.\mc{E}''\right\vert_{D''}.
\end{displaymath} 
One can define a holomorphic bundle over $M_S$ by an induced bundle gluing 
\begin{displaymath}
        \mc{E}^g_S\doteq \mc{E}'\tsum^g \mc{E}''\rightarrow M_S
\end{displaymath}
of the following form. First, fix holomorphic trivialisations over neighbourhoods of infinity $U^{(i)}\subset W^{(i)}_\infty$ along the ends. Then, spreading $g$ along $U^{(i)}$ via pull-back by the fibration maps $\tau^{(i)}:W^{(i)}\rightarrow D^{(i)}$, we identify the fibres of $p_1^*\mc{E}'$ and $p^*_1\mc{E}''$ across the gluing zone of operation $\tsum$. Having said that, I will omit henceforth the superscript $g$ as well as any reference to the particular choice of trivialisations. This paper proves the following \emph{Gluing theorem}:
\begin{theorem}         \label{thm gluing}
Let $\mathcal{E}^{(i)}\rightarrow W^{(i)}$, $i=1,2$, be asymptotically stable bundles of same semi-simple structure group $G$, with \emph{asymptotically rigid} reference metrics $H_0^{(i)}$, admitting a $G-$isomorphism $g:\left.\mc{E}'\right\vert_{D_J'}\tilde\rightarrow \left.\mc{E}''\right\vert_{D''}$. 

There exists $S_0>0$ such that the bundle $\mathcal{E}_S\rightarrow M_S\doteq W'\tsum W''$ admits a $G_2-$instanton, for every $S\geq S_0$.
\end{theorem}

\newpage
\section{Approximate instantons over $M_S$}

In order to produce a $G_2-$instanton over the compact 7-manifold $M_S$, I use cut-offs to obtain a connection on $\mc{E}_S$ which is `approximately' an instanton, then show that a certain surjectivity requirement is satisfied, in order to perturb it into a true solution. 

\subsection{Preliminary moduli theory} 
\label{subsec prelim moduli theory}

Following \cite[pp.83-87]{floer}, given an asymptotically stable bundle ${\mc{E}}\rightarrow W$, we are interested in gauge classes of connections on $\tilde{\mc{E}}=p_1^*\mc{E}\rightarrow M=W\times S^1$ which are `asymptotically HYM'  [cf. (\ref{eq ell complex over D})]. For every (smooth) reference metric $H_0$, we denote its Chern connection $A_0$ and pose [cf. \emph{Notation \ref{Not finite cylinder}}]
\begin{equation}        \label{eq def A}
        \mc{A}\doteq\left\{ p_1^*(A_0+a)
        \left\vert 
        \begin{array}{ll}
        \left. A_0\right\vert_D\in\mc{M}_D, &  \\
        \left\vert a\right\vert, \left\vert \nabla_{A_0} a \right\vert\in L^p_k, &
        a\tublim0 \\
        \end{array}
        \right.\right\},
\end{equation}
taking,
for suitable integers $p,k$ [cf. $(\ref{eq suitable p,k})$ below], the $L_{k}^{p}-$norm induced by $\varphi $:
\begin{equation}       \index{norm!of bundle-valued form}
                        \label{rem Sobolev norm of a}
        \left\Vert f\right\Vert\doteq\left\Vert f\right\Vert _{L_{k}^{p}}
        =\left(\int_{M}\sum_{l=0}^{k}\left\vert \nabla ^{i_{1}}...
        \nabla^{i_{l}}f\right\vert ^{p}d\func{Vol} \right) ^{\frac{1}{p}}.
\end{equation}%
This has in view the use of Sobolev's embedding  (\emph{%
Lemma \ref{Lemma Sobolev's embedding}}) in \emph{Subsection \ref{subsect
subbundle of kernels}}.
For notational clarity I will henceforth leave implicit the pull-back $p_1^*$. 

Posing the gauge-equivalence condition 
\begin{displaymath}
        A_1\sim A_2 \ \Leftrightarrow \ A_2=g(A_1)\overset{loc}{=}A_1-d_1 g.g^{-1}, \quad g\in L^p_{k+1,loc}(\func{Aut}\mc{E}),
\end{displaymath}
and adopt accordingly the gauge group
\begin{displaymath}
        \mc{G}\doteq \left\{ g\in\func{Aut}\tilde{\mc{E}} 
        \left\vert \left\vert\nabla_0g.g^{-1} \right\vert\in L^p_k, 
        \quad g\tublim 1\right.\right\},
\end{displaymath}
whose (bundle of) Lie algebra(s) we denote $\mathfrak{g}$. Since every $\left.A_0\right\vert_D$ is assumed \emph{irreducible} [cf. p.\pageref{def bundle E->W}],  
$\mc{G}$ is in fact a Banach Lie group and the action $\mc{G}\times\mc{A}\rightarrow\mc{A}$ is smooth.
Finally, the Coulomb gauge condition provides transversal slices for the action [cf. (\ref{eq Hodge theory}), below]: 
\begin{equation}     \label{nbhds Te(A)}
U_{\varepsilon }\left( A\right) 
        =A+\left\{ a\in \Omega ^{1}\left( \mathfrak{g}\right) 
        \left\vert \ d_{A}^{\ast }a=0,
        \quad\left\Vert a\right\Vert <\varepsilon, 
        \quad a\tublim0 \right.\right\},
\end{equation}
so that the quotient $\mc{B}=\mc{A}/\mc{G}$ is a Banach manifold \cite[Prop.
4.2.9, p.132]{4-manifolds}:

\begin{proposition}     \label{local model for B}
                        \index{moduli space!local model}
If $A$ is irreducible then, for small $\varepsilon $, the projection from $%
\mathcal{A}$ to $\mathcal{B}$ induces a homeomorphism from  $T_{\varepsilon }\left( A\right)$ to a neighbourhood of $\left[ A\right] $ in $\mathcal{B}$.
\end{proposition}
\begin{definition}      \label{def Moduli space of G2-instantons}
                        \index{moduli space!of G2-instantons@of $G_{2}-$instantons}
The \emph{moduli space of (irreducible) }$G_{2}-$\emph{instantons} on $\mc{E}$ is:%
\begin{equation*}
        \mc{M}^+\doteq \left\{ \left[ A\right]\in \mathcal{B}
        \mid  p_{+}\left( F_{A} \right)\doteq F_A\wedge\ast\varphi=0\right\} .
\end{equation*}
\end{definition}
\noindent NB.: In particular, the instantons obtained in the first paper \cite{G2I}, as solutions of the HYM problem,
decay exponentially in all derivatives to $H_0$ along the cylindrical end, hence $\mc{M}^+\neq \emptyset$ [cf. (\ref{eq HYM solution}) and (\ref{eq def A})]. 

For the local description of $\mc{M}^+$, define around a solution $A$
the map 
\begin{eqnarray}
\begin{array}{r c l}
        \psi \; : \; U_{\varepsilon }\left( A\right) \subset \mathcal{A}
        &\rightarrow&
\Omega^{6}\left( \mathfrak{g}_{E}\right)  \label{psi} \\
a &\mapsto &\psi \left( a\right) \doteq p_{+} \left(F_{A+a}\right)=\left(d_{A}a + a\wedge a\right)\wedge\ast\varphi
\end{array}
\end{eqnarray}%
and write $Z\left( \psi \right) \subset $ $T_{\varepsilon }\left( A\right) $ for its zero set. Slicing out by Coulomb gauge indeed makes $\psi$ a Fredholm map, and we have: 
\begin{proposition}     \index{moduli space!local model}
                        \label{Prop local description of ME}
If $A\in \mathcal{A}$ is an irreducible $G_2-$instanton, then an
$\varepsilon-$neighbourhood of $\left[ A\right] \in \mathcal{M}^+$ is modelled on $Z(\nu)$, where
$\nu$ is the invertible map between finite-dimensional spaces defined by
\begin{gather*}
\nu \; : \;
\underset{
\begin{tabular}{c}
$\cap $ \\ 
$\Omega ^{1}\left( \mathfrak{g}\right) $%
\end{tabular}%
}{\underbrace{\ker \left( d_{A}^{\ast }\oplus d_{A}^{+}\right) }}\;
\longrightarrow \; 
\underset{%
\begin{tabular}{c}
$\cap $ \\ 
\;$ \Omega^{6}\left( \mathfrak{g}\right)$%
\end{tabular}%
}
{\underbrace{\func{coker}d_{A}^{+} \;\cap\; \ker d_{A}}} \\
\quad
\nu \left( a\right) =\sigma \left( 0,a\right)
\end{gather*}%
and $\sigma $ is the non-linear part of the local Fredholm decomposition of $\psi $.
\end{proposition}

The proof of \emph{Proposition \ref{Prop local description of ME}} is postponed to \emph{Section \ref{sect local model moduli space}}, as part of a more detailed discussion of the moduli theory. 

\subsection{Truncating instantons with decaying error term}

We may now start in earnest the proof of \emph{Theorem \ref{thm gluing}}. Let $A'$ be a $G_2-$instanton on $p^*_1\mc{E}'$ [cf. (\ref{eq HYM solution})]. Along a neighbourhood of infinity down the tubular end of $W'$, we write $A'=A_{0}'+ a'$, where $A_0'$ is the lifted Chern connection associated to the reference metric $H_0'$ and  $a'\tublim 0$. Fix a smooth cut-off  
\begin{displaymath}
\begin{array}{c}
        \chi:\R^+\rightarrow \left[ 0,1 \right] \\
        \chi(s)=
        
        \left\{
        \begin{array}{c}
                1, \quad s<0 \\
                0, \quad s\geq 2 \\
        \end{array} .
        \right.\\
\end{array}
\end{displaymath}
and truncate to
\begin{equation*}
        A'_S\doteq A_{0}' +\chi\left( s-S+3 \right)a',
\end{equation*}
which agrees with $A_{0}'$ over the gluing region $[S-1,S]$ and has self-dual part $F_{A'_S}^+\simeq F_{A'_S}\wedge\ast\varphi$ supported in the (topological) cylinder segment $\Sigma_1(S-1)$. Clearly
\begin{displaymath}
        \left\Vert F_{A'_S}\wedge\ast\varphi\right\Vert_{L^p_k\left( M'\right)}
        \leq C_{p,k}\e{-S}
\end{displaymath}
for $M'\doteq W'\times S^1$.
Repeating the construction for $A''_S$ on $p_1^*\mc{E}''$, we may assume that the truncated connections $A^{(i)}_S$ match (via the bundle isomorphism $g$) over $[S-1,\infty[$, hence in particular over $[S-1,S]$, so they glue together to define a smooth connection
\begin{displaymath}
        A_S\doteq A'\tsum A'' \quad \text{on} \quad \mc{E}_S.
\end{displaymath}
Since the hyper-K\"ahler rotation at infinity is assumed to be an isometry, we still have the asymptotic decay of the self-dual part
\begin{displaymath}
        \left\Vert F_{A_S}\wedge\ast\varphi \right\Vert_{L^p_k\left(M_S\right)}
        \leq C_{p,k}\e{-S}
\end{displaymath}
from which we see that $A_S$ is \emph{almost} an instanton. 

\subsection{Non-degeneracy under acyclic limits}

Let $A=A^{(i)}$ as above, say, be a $G_2-$instanton on
the pull-back bundle\begin{displaymath}
        \tilde{\mc{E}}=p_1^*\mathcal{E}\rightarrow M
        \doteq W\times S^1.
\end{displaymath}
Then the hypothesis that the connection at infinity is \emph{acyclic} implies that $A$ itself is acyclic, with respect to its own deformation complex:  
\begin{equation}        \label{extended complex 2}
        \Omega ^{0}\left( \mathfrak{g}\right) \;
        \overset{d_{A}}
        {\underset{d_{A}^{\ast }}{\rightleftarrows }}\;
        \Omega ^{1}\left( \mathfrak{g}\right) \;
        \overset{d_{A}^{+}}
        {\overbrace{\overset{d_{A}}{\longrightarrow }\;
        \Omega^{2}\left( \mathfrak{g}\right)  \;
        \overset{\ast \varphi \wedge .}
        {\longrightarrow }}}\;
        \Omega ^{6}\left( \mathfrak{g}\right) \;
        \overset{d_{A}}
        {\underset{d_{A}^{\ast }}{\rightleftarrows }}\;
        \Omega ^{7}\left( \mathfrak{g}\right).
\end{equation}%
The goal of this \emph{Subsection} is to prove that claim, in the following terms:

\begin{proposition}     \label{Prop asymp acyclic => acyclic}
When the reference metric $H_0$ is asymptotically rigid, then the induced $G_2-$instanton $A$ lifted from a HYM solution [cf. $(\ref{sd part kahler -> G2})$ and $(\ref{eq HYM solution})$] is  acyclic, i.e., $\mathbf{}\mathbf{H}_A^0=0$ and $\mathbf{H}_A^1=0$ in the deformation complex $(\ref{extended complex 2})$.
\end{proposition}
Since our asymptotically stable bundle is, by definition, indecomposable, we have already $\mathbf{H}_A^0=0$, so one only needs to check the non-degeneracy condition $\mathbf{H}_A^1=0$. On the other hand, the complex is self-dual (under the Hodge star), so this is equivalent to showing $\mathbf{H}_A^2=0$, which means precisely that $A$ is an isolated point in its moduli space $\mathcal{M}^+$, in the light of the local model given by \emph{Proposition $\ref{Prop local description of ME}$}.
To check this fact, we resort to the Chern-Simons functional
\begin{displaymath}
        \rho(b)_A=\int_{W\times S^1}\tr F_A\wedge b_A\wedge \ast\varphi.
\end{displaymath} 

For  any given direction $a\in T_{[A]}\mathcal{B}\;\widetilde{\subset}\;\Omega^1(\mathfrak{g}) $ with $\Vert a \Vert=1$ and possibly short length $\varepsilon>0$, we have
\begin{eqnarray*}
        \left[A+ha\right]\in\mc{M}^+&\Leftrightarrow& 
        0\overset{\forall b}{=}\rho(b)_{A+\varepsilon a}=\underset{0}{\underbrace{\rho(b)_A}}
        +\varepsilon.D\left[ \rho(b) \right]_A(a)+O(\varepsilon^2)
\end{eqnarray*} 
for vector fields $b\in \Gamma\left( T\mathcal{A} \right)$. Explicitly, the first order variation is
\begin{eqnarray*}
        D\left[ \rho(b) \right]_A(a)&=&\int_{W\times S^1}\tr\{d_Aa\wedge b_A 
        + (Db)_A(a)\wedge 
        \underset{F_A^+=0}{\underbrace{F_A\}\wedge\ast\varphi}}\\
        &=&\int_{W\times S^1}\tr a\wedge d_A b_A\wedge\ast\varphi 
        +\underset{0}{\underbrace{\lim_{S\rightarrow\infty}\int_{\partial W_S\times S^1}
        \tr \left\{ a\wedge b_A \right\}\wedge\ast\varphi}}\\
        &=&\int_{W\times S^1}\tr a\wedge d_A^+b_A
\end{eqnarray*}
since $H_0$ is asymptotically rigid and so $\left\vert a\right\vert\underset{S\rightarrow\infty}{\longrightarrow}0$, for small enough $\varepsilon$. Notice in passing that this is zero for any direction $a\in\img d_A\subset\Omega^1(\mathfrak{g})$ along gauge orbits, corresponding to the intuitive fact that the only `meaningful' perturbations are those which descend nontrivially to $\mc{B}$. Choose now $0\neq\xi\in\Omega^6(\mathfrak{g})$ such that 
\begin{displaymath}
        b_A\doteq\left(d_A^+\right)^*\xi\neq0 
        \text{ in } 
        \Omega^1(\mathfrak{g}).
\end{displaymath}
By the orthogonal decomposition $(\ref{eq Hodge theory})$, this can be done in such a way that $d_A^+b_A\neq0$, so for any direction $a\in\Omega^1(\mathfrak{g})$ (transverse to gauge orbits) the number
\begin{displaymath}
        D\left[ \rho(b) \right]_A(a)\doteq N_A(a,\xi)\in\R
\end{displaymath}
is not zero for a generic choice of $\xi$. Rescaling $\widetilde{b_A}\doteq \frac{1}{N_A(a,\xi)}b_A$, we find 
\begin{displaymath}
        \rho(\tilde{b})_{A+ha}
        =\varepsilon+O(\varepsilon^2)O(\Vert\tilde b_A \Vert)\neq0,
        \quad \text{for } \;0\neq \varepsilon\ll 1.
\end{displaymath}

Now, since $\mathcal{M}^+$ is finite-dimensional, there are tangent vectors ($1-$forms on the base)  $u_1,\dots,u_n\in T_{[A]}\mathcal{B}$ such that any such perturbation is written as
\begin{displaymath}
        \varepsilon.a=\varepsilon.(a^1u_1+\dots+a^n u_n), \quad a^i\in \R.
\end{displaymath}
But in the above way we can find, respectively for $u_1,\dots,u_n$, vector fields $\tilde b_1,\dots,\tilde b_n$ such that $\rho(\tilde{b_i})_{A+hu_i}= \varepsilon+O(\varepsilon^2)$. Consequently, for a generic linear combination $\mathbf{\tilde b}\doteq \beta^1 \tilde b_1+\dots+\beta^n \tilde b_n$, one has
\begin{displaymath}
\rho(\mathbf{\tilde{b}})_{A+ha}=\varepsilon.\underset{\neq0}{\underbrace{\left( a^1\beta^1+...+a^n\beta^n \right)}}+O(\varepsilon^2).
\end{displaymath}Hence there exists a (possibly small) value  $\varepsilon_{0}>0$ such that $\rho_{A+\varepsilon a}\neq0$, as a $1-$form on $\mathcal{A}$, for any $\varepsilon.a\in U_{\varepsilon_0}(\left[A\right])$. In other words, there are no instantons in the open ball of radius $\varepsilon_{0}$ around $\left[A\right]$ in the moduli space (\emph{q.e.d.}). 

\newpage
\section{Perturbation theory over long tubular ends}

Following the standard approach, we may now look for a nearby exact solution $A=A_S+a$ to the $G_2-$instanton equation 
\begin{equation}        \label{eq nearby instanton}
        d_{A_S}^+a+\left( a\wedge a \right)\wedge\ast\varphi
        =-F_{A_S}\wedge\ast\varphi\doteq\epsilon(S).
\end{equation}
We adopt all along the \emph{acyclic assumption}  that the operators $d_{A\ii}^+$ have trivial cokernel, i.e., the (irreducible) connections $A^{(i)}$ are isolated points in the respective moduli spaces $\mathcal{M}^+_{\mathcal{E}^{(i)}}$ 
[cf. \emph{Proposition \ref{Prop asymp acyclic => acyclic}}].

\subsection{Noncompact Sobolev estimates}

For briefness, let us refer to ordered integers $k> l\geq0$ and $q\geq p$ as \emph{suitable} if they satisfy the Sobolev condition:
\begin{equation}        \label{eq suitable p,k}
        \frac{1}{p}-\frac{1}{q}\leq \frac{k-l}{7}.
\end{equation} 
In particular, the pair  $p,k$ will be \emph{suitable} when  $k=l+1$ and $2p=q$ are suitable.

\begin{lemma}   \label{Lemma noncompact Sobolev}
Let $ W$ be an asymptotically cylindrical $3-$fold; given suitable $k\geq l$ and $q\geq p$, there exists a constant $C\doteq C_{W,p,q,k,l}>0$ such that, for sections of any bundle over $W\times S^1$ (with metric and compatible connection), 
\begin{displaymath}
        \left\Vert .\right\Vert_{L^q_l}\leq 
        C\left\Vert . \right\Vert_{L^p_k} .
\end{displaymath}
\begin{proof}
Following \cite[pp.70-72]{floer}, set $B_0\doteq W_0\times S^1$ and consider for $n\geq1$ the tubular segments of `length one' $B_n=\left(W_n\setminus W_{n-1}\right)\times S^1$ along the tubular end. Then, using the usual Sobolev estimate for compact domains,
\begin{eqnarray*}
        \left\Vert f \right\Vert_{L^q}^q &=& 
        \sum_{n\in \mathbb{N}} \int_{B_n} \left\vert f \right\vert^q\\
        &\leq&\sum_{n\in \mathbb{N}}  C_{n}\left( \int_{B_n}
        \left\vert \nabla f \right\vert^p 
        + \left\vert f \right\vert^p\right)^{q/p}\\
        &\leq& \tilde C \left( \sum\int\left\vert \nabla f \right\vert^p
        + \left\vert f \right\vert^p \right)^{q/p} 
        = \tilde C \left\Vert f \right\Vert^q_{L^p_1}
\end{eqnarray*}
with $\tilde C=\lim\sup C_n<\infty$, since the segments are asymptotically cylindrical. This proves the statement, by induction on $l$. 
\end{proof}
\end{lemma}

\newpage
\subsection{Gluing right inverses}

We now  investigate the behaviour of right-inverses under truncation and gluing: 

\begin{lemma}
For $S\gg0$, the operators $d_{A_S\ii}$ admit bounded right inverses $Q_S\ii$ satisfying
\begin{displaymath}
        \Vert Q_S\ii\xi\Vert_{L^p_{k}}
        \leq C\ii_{p,k}\left\Vert \xi \right\Vert_{L^p_{k-1}}
\end{displaymath}
for suitable $p,k\in\mathbb{N}$, where the bound $C\ii_{p,k}$ depends only on $A\ii$, not on $S$.
\begin{proof}
The operators  $d_{A\ii}$ correspond to the original instantons over each tubular component $W\ii\times S^1$, hence by the acyclic assumption  [cf. \emph{Proposition \ref{Prop asymp acyclic => acyclic}}] they admit bounded right inverses $Q\ii$,  independent of $S$. 

The crucial fact is that $a^{i}_S\doteq A_S\ii-A\ii=O\left( \e{-S} \right)$ so, for $S\gg0$, a right inverse for $d_{A\ii}$ gives a right inverse $Q_S\ii$ for $d_{A\ii_S}$ with (approximately) the same uniform Lipschitz bound.
\end{proof}
\end{lemma}

\begin{corollary}
For $S\gg0$, there exist an `approximate' right inverse $Q_{S}$ and a true right inverse $P_{S}$ for $d_{A_S}^+$:
\begin{displaymath}
        P_{S}=Q_{S}\left( d_{A_S}^+Q_{S} \right)^{-1}.
\end{displaymath}
Moreover, for suitable $p,k\in\mathbb{N}$, there is a uniform bound $C_{p,k}$ on the operator norm of $P_{S}$:
\begin{displaymath}
        \left\Vert P_{S}(\xi)\right\Vert_{L^p_{k}}\leq C_{p,k}
        \left\Vert \xi \right\Vert_{L^p_{k-1}}.
\end{displaymath}
\begin{proof}
Following a standard argument \cite[\S3.3 \& \S4.4]{floer}, we may take truncation functions $\chi\ii_S:W\ii\rightarrow \left[ 0,1 \right]$ satisfying
\begin{displaymath}
        \left(\chi_{S}'\right)^2+\left( \chi_{S}'' \right)^2=1,
        \quad \func{supp}\chi_{S}\ii\subset W\ii_{\frac{3S}{2}},         
        \quad \Vert\nabla\chi_{S}\ii\Vert_{L^\infty}=O(\e{-S}).
\end{displaymath} 
Then, denoting $r_{S}\ii:M_S\rightarrow W\ii\times S^1$ the maps  given by restriction for $s\leq2S$ and extended by zero along the rest of the tubular end, we form
\begin{equation}
        Q_{S}\doteq\left(\chi'_{S}\right)^2 \left(Q'_S\circ r'_S\right)
        +\left(\chi''_{S}\right)^2 \left(Q''_{S}\circ r''_S\right)
\end{equation}
and we check that this is an approximate right inverse in the sense that
$\Vert d_{A_S}^+Q_S-I\Vert=O(\e{-S})$; indeed we have
\begin{displaymath}
        d_{A_S}^+Q_S =\underset{I}{\underbrace{ 
        \sum\nolimits_{i=1,2} \{\left(\chi\ii_S\right)^2d_{A_S}^+ \left(Q\ii_{S}\circ r\ii_S \right)}
        }+ 2\underset{O(\e{-S})}{\underbrace{ \left( \chi\ii_S\nabla\chi\ii_S\right)
        \diamondsuit\left( Q\ii_{S}\circ r\ii_S \right)}}\}
\end{displaymath}
where $\diamondsuit$ denotes an algebraic operation. It follows from the  \emph{Lemma}  that the second summand is dominated by the decay of $\Vert\nabla\chi_{S}\ii\Vert_{L^\infty}$.
Then
\begin{displaymath}
         P_S\doteq Q_S \left( d_{A_S}^+Q_S \right)^{-1}
\end{displaymath}
is a true right inverse for $d_{A_S}^+$, with uniformly bounded norm determined by $Q_S$ and $d_{A_S^+}$ itself.
\end{proof}
\end{corollary}

\subsection{Exact solution via the contraction principle}
For a solution of the form $a=P\xi$, equation $(\ref{eq nearby instanton})$ reads
\begin{equation}        \label{eq contraction}
        \left(I+G\right)(\xi)=\epsilon(S)
\end{equation}
where $G(\xi)\doteq  P\left( \xi\right)\wedge P\left( \xi\right)
        \wedge\ast\varphi$ and $\left\Vert\epsilon(S)\right\Vert$ is small, as $S\gg0$.
Thus if, given  suitable $p,k$, the map [cf. \emph{Lemma \ref{Lemma noncompact Sobolev}} for last inclusion]
\begin{displaymath}
        G: L^{2p}_{k-1}\left(\Omega^6\right)
        \hookrightarrow L^{p}_{k-1}\left(\Omega^6\right)
        \longrightarrow L^{p}_{k}\left(\Omega^6\right)
        \hookrightarrow L^{2p}_{k-1}\left(\Omega^6\right)
\end{displaymath}
is Lipschitz with constant strictly smaller than $1$, then by the Banach contraction principle $I+G$ is continuously invertible between neighbourhoods of $\left(I+G\right)(0)=0$ and we obtain, for large $S$ say, a solution as
\begin{displaymath}
        a=P\left( I+G \right)^{-1}\epsilon(S).
\end{displaymath} 

In order to prove this, we need a uniform bound on exterior multiplication, which is an immediate consequence of H\"older's inequality:

\begin{proposition}     \label{prop bounded multiplication}
For suitable $p,k\in\mathbb{N}$, there exists a uniform constant $M_{p,k}$ such that 
\begin{displaymath}
        \left\Vert a\wedge b \right\Vert_{L^{p}_{k-1}}\leq M_{p,k} 
        \left\Vert a \right\Vert_{L^{2p}_{k-1}}
        \left\Vert b \right\Vert_{L^{2p}_{k-1}}.
\end{displaymath}
\end{proposition}

From this we infer that operator $G$ is bounded uniformly in $S$ for any suitable Sobolev norm, but so far we have no definite control over the actual bound. We can scale away this apparent difficulty using the fact that $G$ is homogeneous of degree 2; letting $\lambda=MC^{2}$, $\tilde\xi=\lambda\xi$ and $\tilde G=\frac{1}{\lambda}G$, equation $(\ref{eq contraction})$ becomes
\begin{displaymath}
        \left( I+\tilde G \right)(\tilde\xi)=\lambda\epsilon(S).
\end{displaymath}
The point here is that now $\tilde G$ is a contraction over the ball $\tilde B_{1}\doteq\left\{ \Vert \tilde\xi\Vert\leq1 \right\}$:
\begin{displaymath}
        \Vert \tilde G (\tilde\xi) \Vert \leq
        \frac{MC^{2}}{\lambda}\Vert \tilde\xi \Vert^2
        \leq \Vert \tilde\xi\Vert.
\end{displaymath}
Then indeed $I+\tilde G$ is a homeomorphism onto  an interior domain $0\in U\subset \tilde B_1$, and one can choose $S\gg0$ so that $\lambda\epsilon(S)\in U$. One may check, by bootstrapping \cite[.96]{floer}, that  $A$ is in fact smooth. We have thus proved \emph{Theorem \ref{thm gluing}}.

In conclusion, it should be noticed that the acyclic hypothesis is a rather strong, non-generic requirement. The fact that it is not void should therefore be illustrated:

\begin{example}

The prime Fano $3-$folds of type $X_{22}$ were discovered by Iskovskikh \cite{X22} and further studied by Mukai \cite[\S3]{Mukai Fano 3-folds}. On one hand, these appear in Kovalev's list of suitable blocks for the gluing construction; in fact, they are extremal in the sense that a pair of base manifolds of type $X_{22}$ realises the lower bound on the third Betti number $b_3(M)=71$ \cite[pp.158-159]{kovalevzao}.

On the other hand, crucially, these come equipped with an asymptotically stable bundle $E\rightarrow X_{22}$ which is \emph{rigid} \cite[\S3]{Mukai Moduli of bundles} over a divisor $D\in\left\vert -K_{X_{22}} \right\vert$. In other words, the holomorphic bundle $\left.E\right\vert_D$  corresponds to an isolated point in its moduli space, hence its associated HYM metric $H_0$ is indeed \emph{acyclic}.
\end{example}

\section{Local model for the moduli space}
\label{sect local model moduli space}

The moduli space $\mathcal{M}^+$ of $G_2-$instantons on $\mc{E}\rightarrow M$  is locally described as the zero set of a map $\psi $ [cf. $\left( \ref{psi}\right) $] between the Banach
spaces $U_{\varepsilon }\left( A\right) \subset \mathcal{A}$ and 
$\Omega^6\left( \mathfrak{g}\right) $. Therefore, if our map $\psi $ is Fredholm on $Z\left( \psi \right) $, it is a matter of standard theory to model a
neighbourhood of $\left[ A\right] $ in $\mathcal{M}^+$ on the
finite-dimensional set $\nu ^{-1}\left( 0\right) $ [cf. 
\emph{Corollary \ref{cor zero set}} in the \emph{Appendix}] . Foreseeing \emph{Proposition \ref{prop projection between kernels}}, we 
restrict attention to  \emph{irreducible} connections on an $SU\left(n\right) -$bundle  $\mc{E}$.

\subsection{Noncompact Fredholm theory}
\label{Subsec Fredholm theory}
\label{subsect extended elliptic complex}
\label{subsect subbundle of kernels}

As defined before, the map $\psi $
is just the self-dual part of the curvature, so $\psi \left( a\right) -\psi
\left( 0\right) =\left(p_+\circ d_{A}\right)a + O( \left\vert a\right\vert ^{2}) $ and%
\begin{equation*}
        \left( D\psi \right) _{0}=p_+\circ d_{A}.
\end{equation*}
Moreover, by the `slicing' condition $\left( \ref{nbhds Te(A)}%
\right) $ across orbits, we consider in fact the restriction%
\begin{equation}        \label{da+ restr kerda*}
                        \index{gauge orbit}
        \begin{array}{ccccc}
        p_+\circ d_{A} & : & \ker d_{A}^{\ast } & \longrightarrow 
        & \Omega _{+}^{2}\left( \mathfrak{g}\right) . \\
        &  & \cap\ &  &  \\
        &  & \Omega ^{1}\left( \mathfrak{g}\right)  &  &  \\
\end{array}
\end{equation}

Since the map $L_{\ast \varphi }=`\ast \varphi \wedge.$' now plays the role of `SD projection',  we denote henceforth 
\begin{equation}        \label{Da+}
        d_{A}^{+}=L_{\ast \varphi }\circ d_{A}:\Omega ^{1}
        \left( \mathfrak{g}\right) \rightarrow \Omega ^{6}
        \left( \mathfrak{g}\right)
\end{equation}
and consider the extended deformation complex%
\index{elliptic!complex}%
\begin{equation}        \label{extended complex}
        \Omega ^{0}\left( \mathfrak{g}\right) \;
        \overset{d_{A}}
        {\underset{d_{A}^{\ast }}{\rightleftarrows }}\;
        \Omega ^{1}\left( \mathfrak{g}\right) \;
        \overset{d_{A}^{+}}
        {\overbrace{\overset{d_{A}}{\longrightarrow }\;
        \Omega^{2}\left( \mathfrak{g}\right)  \;
        \overset{\ast \varphi \wedge .}
        {\longrightarrow }}}\;
        \Omega ^{6}\left( \mathfrak{g}\right) \;
        \overset{d_{A}}
        {\underset{d_{A}^{\ast }}{\rightleftarrows }}\;
        \Omega ^{7}\left( \mathfrak{g}\right).
\end{equation}
Using $d\ast \varphi =0$, we find
\begin{equation}
        \left[ L_{\ast \varphi },d_{A}\right] =0,  \label{[Lfi,dA]=0}
\end{equation}%
so, when $A$ is an instanton, $\left( \ref{extended complex}\right) $ is
indeed a complex and the identification of the self-dual $2-$forms with the $6-$forms is consistent with the relevant differential operators (for more on elliptic complexes under the condition $d*\varphi=0$, see \cite{Fernandez&Ugarte}). Moreover, this complex is elliptic:\begin{lemma}
\label{lemma (Da+)*=*(Da+)*}The operator $d_{A}^{+}$ defined by $\left( \ref%
{Da+}\right) $ has formal adjoint 
\begin{equation*}
\left( d_{A}^{+}\right) ^{\ast }
=\ast d_{A}^{+}\ast :\Omega ^{6}\left(\mathfrak{g}\right) 
\rightarrow \Omega ^{1}\left( \mathfrak{g}\right) .
\end{equation*}
\begin{proof}
For $a\in \Omega ^{1}\left( \mathfrak{g}\right) $ and $\eta \in \Omega
^{6}\left( \mathfrak{g}\right) $,
we have pointwise:
\begin{eqnarray*}
\left\langle d_{A}^{+}a,\eta \right\rangle \left( \ast 1\right) &=&\left(
\ast \varphi \wedge d_{A}a\right) \wedge \ast \eta =\left( d_{A}a\right)
\wedge \ast \left( \ast \left( \ast \varphi \wedge \ast \eta \right) \right)
\\
&=&\left\langle d_{A}a,\ast \left( L_{\ast \varphi }\ast \eta \right)
\right\rangle =\left\langle a,d_{A}^{\ast }\left( \ast L_{\ast \varphi }\ast
\eta \right) \right\rangle \left( \ast 1\right) \\
&=&\left\langle a,\ast \left( d_{A}L_{\ast \varphi }\right) \ast \eta
\right\rangle \left( \ast 1\right) \overset{\left( \ref{[Lfi,dA]=0}\right) }{%
=}\left\langle a,\ast \left( L_{\ast \varphi }d_{A}\right) \ast \eta
\right\rangle \left( \ast 1\right) \\
&=&\left\langle a,\left( \ast d_{A}^{+}\ast \right) \eta \right\rangle
\left( \ast 1\right) 
\end{eqnarray*}
\vspace{-1pt}
\end{proof}
\end{lemma}

\begin{proposition}     \label{prop elliptic complex}
When $A$ is a $G_2-$instanton, the complex $\left(\ref{extended complex}\right)$  is elliptic.
\begin{proof}
First of all, since $\left( d_{A}^{+}\right) ^{\ast }=\ast d_{A}^{+}\ast $ [%
\emph{Lemma \ref{lemma (Da+)*=*(Da+)*}}] and $d_{A}^{\ast }=\ast d_{A}\ast $%
, notice that our complex is self-dual with respect to the Hodge star:%
\begin{equation*}
\begin{array}{ccccccc}
        \Omega ^{0} & \overset{d_{A}}{\longrightarrow } & \Omega ^{1} & 
        \overset{d_{A}^{+}}{\longrightarrow } & \Omega ^{6} & 
        \overset{d_{A}}{\longrightarrow } & \Omega ^{7} \\ 
        \shortparallel & & \shortparallel &  & \shortparallel & & \shortparallel\\
        \ast \Omega ^{7} & \underset{d_{A}^{\ast }}{\longleftarrow } & 
        \ast \Omega ^{6} & \underset{\left( d_{A}^{+}\right) ^{\ast }}
        {\longleftarrow } & \ast\Omega ^{1} & \underset{d_{A}^{\ast }}
        {\longleftarrow } & \ast \Omega ^{0}%
\end{array}%
.
\end{equation*}%
By \emph{Corollary \ref{cor complex is elliptic iff dual is elliptic}}, it
suffices to show ellipticity at $\Omega ^{1}\left( \mathfrak{g}\right) $%
, as that is equivalent to the ellipticity of the dual $\ast \Omega ^{7}%
\overset{d_{A}^{\ast }}{\longleftarrow }\ast \Omega ^{6}\overset{\left(
d_{A}^{+}\right) ^{\ast }}{\longleftarrow }\ast \Omega ^{1}$, which is just $%
\Omega ^{1}\overset{d_{A}^{+}}{\longrightarrow }\Omega ^{6}\overset{d_{A}}{%
\longrightarrow }\Omega ^{7}$. Fixing a section $\xi $ of $T^{^{\prime }}M$
(the cotangent bundle minus its zero section), we have symbol maps%
\begin{equation*}
        0\rightarrow \pi ^{\ast }
        \left( \Omega ^{0}\left( \mathfrak{g}\right)\right) _{\xi }
        \overset{\xi .\left( .\right) }{\longrightarrow }\pi ^{\ast }
        \left( \Omega ^{1}\left( \mathfrak{g}\right) \right) _{\xi }
        \overset{\ast \varphi \wedge \xi \wedge \left( .\right) }
        {\longrightarrow }\pi ^{\ast }\left( \Omega ^{6}
        \left( \mathfrak{g}\right) \right) _{\xi}\longrightarrow \dots
\end{equation*}%
For $\alpha \in \Omega ^{1}\left( \mathfrak{g}\right) $ such that $%
\ast \varphi \wedge \xi \wedge \alpha =0$, exactness means $\alpha $ has to
lie in $\xi .\Omega ^{0}\left( \mathfrak{g}\right) $. Since $G_{2}$ acts
transitively on $S^{6}$, take $%
g\in G_{2}$ such that $g^{\ast }\xi =\left\Vert \xi \right\Vert .e^{1}$ and
denote $\widetilde{\alpha }=g^{\ast }\alpha $, so that%
\begin{equation*}
        \ast \varphi \wedge e^{1}\wedge \widetilde{\alpha }=0.
\end{equation*}%
That is just the statement that $e^{1}\wedge \widetilde{\alpha }$ is
anti-self-dual, but this cannot occur unless $e^{1}\wedge \widetilde{\alpha }%
=0$, as%
\begin{equation*}
        \left( e^{1}\wedge \widetilde{\alpha }\right) \wedge \varphi 
        =\widetilde{\alpha }\wedge \left( e^{1567}-e^{1345}-e^{1426}+e^{1237}\right)
\end{equation*}%
has non-vanishing components involving $e^{1}$ and $\ast \left( e^{1}\wedge 
\widetilde{\alpha }\right) $ obviously has not. Therefore $\widetilde{\alpha 
}=f.e^{1}$ for some $f\in \Omega ^{0}\left( \mathfrak{g}\right) $, and%
\begin{equation*}
        \alpha =\left( g^{\ast }\right) ^{-1}\left( f.e^{1}\right) 
        =\frac{f}{\left\Vert \xi \right\Vert }.\xi \in \xi .\Omega ^{0}
        \left( \mathfrak{g}\right) .
\end{equation*}
\end{proof}
\end{proposition}

In view of the isomorphism  $\left. L_{\ast \varphi} \right\vert _{\Omega
_{+}^{2}}:\Omega_{+}^{2} \, \tilde{\rightarrow} \;\Omega^{6}$, taking the self-dual part of curvature via the $\mathcal{G}-$equivariant map 
$
L_{\ast \varphi}\circ
F^{+}:\mathcal{A}\rightarrow \Omega^{6}\left( \mathfrak{g}\right)$
 defines a section $\Psi \left( \left[ A%
\right] \right) =F_{A} \wedge\ast\varphi $ of the Hilbert bundle 
\begin{equation*}
        \mathcal{A}\times _{\mathcal{G}}\Omega ^{6}
        \left( \mathfrak{g}\right) \rightarrow \mathcal{B}.
\end{equation*}

The \emph{intrinsic derivative}, i.e., the component of the total derivative
tangent to the gauge-fixing slices in $\Omega ^{1}\left( \mathfrak{g}_{E}\right)$, of $\Psi $ at $\left[ A\right] $ is%
\begin{eqnarray*}
        \left( D\Psi \right) _{\left[ A\right] } :\ker d_{A}^{\ast }
        &\rightarrow&\ker d_{A}\subset \Omega^{6}\left( \mathfrak{g}\right) \\
        a&\mapsto& d_{A}^{+}a.
\end{eqnarray*}%
To see that $\left( D\Psi \right) _{\left[ A\right] }$ is Fredholm over $Z\left( \Psi \right) $, consider the extended operator
\begin{equation*}       \index{elliptic!operator}
\begin{array}{c c l}
        {\mathbb{D}_{A}}: \;\Omega ^{1}\left( \mathfrak{g}\right)         \oplus \Omega ^{7}\left( \mathfrak{g}\right) 
        &\rightarrow& \Omega ^{0}\left( \mathfrak{g}\right) 
        \oplus \Omega^{6}\left( \mathfrak{g}\right)\\
        \qquad \left( a,f\right) &\mapsto& 
        \left( d_{A}^{\ast}a,\,d_{A}^{+}a+d_{A}^{\ast }f\right) .
\end{array}
\end{equation*}
If the base manifold was compact, then standard elliptic theory would imply ${\mathbb{D}_{A}}=d_{A}^{\ast }\oplus \left( d_{A}^{+}\oplus
d_{A}^{\ast }\right) $ is Fredholm, as the `Euler characteristic' of an elliptic complex, by  \emph{Proposition \ref{prop elliptic complex}}. However, on our manifolds with cylindrical ends $W\simeq W_0\cup W_\infty$,  the parametrix patching method over the compact piece $W_0$ must be combined with tubular theory over $W_\infty$ under the acyclic assumption. This follows in all respects the proof of \cite[Prop. 3.6]{floer} and its preceding discussion, except that in this case we are in the much simpler situation where the exponential decay of $(a,f)\in\ker\mathbb{D}_{A}$ is guaranteed from the outset by our definitions [cf. (\ref{eq def A})] and the fact that the bundle is indecomposable. Then we have:  
\begin{claim}
If $\left[ A\right] \in Z\left( \Psi \right)$ is asymptotically rigid [cf. \emph{Definition \ref{def asymp rigid}}],  ${\mathbb{D}_{A}}$ is a Fredholm operator.
\end{claim}

In particular, ${\mathbb{D}_{A}}$ has closed range, so there is an orthogonal decomposition:
\begin{equation}        \label{eq Hodge theory}
\Omega^6(\mathfrak{g})=\ker d_A\oplus \img d^*_A
\end{equation} 
and this implies that $\left( D\Psi \right) _{\left[ A\right]} $ is also Fredholm:

\begin{displaymath}
\begin{array}{r c c c c c c }
\ker\left(D\Psi\right)_{\left[ A\right]} & 
\hookrightarrow & \ker{\mathbb{D}_{A}} &  &  &  &  \\
\func{coker}\left( D\Psi \right) _{\left[ A\right] } 
& = & \func{coker}d_{A}^{+} & \cap & 
\func{coker}\left({ d_{A}^*}
\vert_{\Omega ^{7} }\right) & 
\hookrightarrow & \func{coker}{\mathbb{D}_{A.}} \\
 &  & \cap &  & \| &  &  \\
 &  & \ker d_{A} & = & \ker d_{A} &  &  \\
\end{array}
\end{displaymath}

Finally, the moduli space
of $G_2-$instantons [\emph{Definition \ref{def Moduli space of G2-instantons}}] is cut out as its zero set $Z(\Psi)$. As an immediate consequence of $\left( \ref{[Lfi,dA]=0}\right) $ and the
Bianchi identity, we have $\Psi \left( \left[ A\right] \right) \in \ker
d_{A}\subset \Omega ^{6}\left( \mathfrak{g}_{E}\right)$,
so, intuitively, the image of $\Psi $ lies in the `subbundle of kernels of $d_{A}$':
\begin{equation}         \label{subbundle of kernels}
        \begin{array}{ccccc}
        \mathcal{A} & \times _{\mathcal{G}} & \ker d_{A} & \rightarrow  &  \mathcal{B} \\
         &  & \cap &  &  \\
         &  & \Omega ^{6}\left( \mathfrak{g}\right) &  & 
\end{array}
\end{equation}When $\mc{E}$ is an $SU\left( n\right)
-$bundle, say, one can use the orthogonal projections $p_{a}:\ker
d_{A}\rightarrow \ker d_{A_{0}}$ to trivialise the fibres onto $\ker
d_{A_{0}}$ over a neighbourhood $U_{\varepsilon }\left( \left[ A_{0}%
\right] \right) \subset \mathcal{B}$, where $\varepsilon $ is a small global
constant:
\newpage
\begin{proposition}     \label{prop projection between kernels}
Let $\mc{E}$ be an $SU\left( n\right) -$bundle over a compact $G_{2}-$manifold $\left( M,\varphi
\right) $ and $A_{0}$ an irreducible connection; then there exists $\varepsilon >0$ such that the orthogonal projection%
\begin{equation*}
p_{a}:\ker d_{A}\rightarrow \ker d_{A_{0}}
\end{equation*}%
in $\Omega ^{6}\left( \mathfrak{g}\right) $ is an isomorphism for
all $A=A_{0}+a\in U_{\varepsilon }\left( A_{0}\right) $.

\begin{proof}
We consider, throughout, the operators [cf. $\left( \ref{extended complex}\right) $ \& $(\ref{eq Hodge theory})$]:
\begin{equation*}
        \Omega^{6}\left( \mathfrak{g}\right)
        \overset{d_{A_{0}},d_{A}}{\underset{d_{A_{0}}^{\ast },
        d_{A}^{\ast }}{\rightleftarrows }}\Omega^{7}
        \left( \mathfrak{g}\right).
\end{equation*}%
Writing  $\rho =d_{A_{0}}^{\ast }f$ for some $f\in \Omega ^{7}\left( \mathfrak{g%
}\right) $, we denote  elements of $\Omega^6 \left( \mathfrak{g}\right) $ by 
\begin{equation*}
        \eta =(\eta _{0}\oplus \rho )\in 
        \left(\ker d_{A_{0}}\oplus \img d_{A_{0}}^{\ast }\right)
        =\Omega^6\left( \mathfrak{g}\right).
\end{equation*}%

\noindent\emph{\textbf{Surjectivity}}

Given $\eta _{0}\in \ker d_{A_{0}}$, write $g_{0}\doteq
-a\wedge \eta _{0}$; surjectivity of $p_{a}$ means finding $\rho \in 
\img d_{A_{0}}^{\ast }\subset \Omega^6\left( \mathfrak{g}%
\right) $ such that $\eta =\eta _{0}\oplus \rho \in \ker d_{A}$, i.e.,
solving for $\rho $ the equation%
\begin{equation}
d_{A}\rho =g_{0}.  \label{surjectivity}
\end{equation}%
Since $A_{0}$ is irreducible, one has $\left( \img d_{A_{0}}\right) ^{\perp
}=\ker d_{A_{0}}^{\ast }=\left\{ 0\right\} $, therefore $\func{img}d_{A_{0}}=\Omega ^{7}\left( \mathfrak{g}\right)$. Thus one may think of the restriction of $d_{A}$ to $\func{img}%
d_{A_{0}}^{\ast }$ as
\begin{equation}        \label{restriction of dA}
        d_{A}:\img d_{A_{0}}^{\ast }\rightarrow \img d_{A_{0}}.
\end{equation}%
Bijectivity of linear maps between Banach spaces is an open condition
[\emph{Lemma \ref{Lemma Fine}}], so one can show that $\left( \ref%
{restriction of dA}\right) $ is invertible by checking that, for suitably
small $a$, this map is arbitrarily close to the isomorphism $d_{A_{0}}:\img d_{A_{0}}^{\ast }\tilde{\rightarrow}\;\img d_{A_{0}}$. Indeed, writing $L_{a}:\eta \mapsto a\wedge \eta$, there exists a global constant such that%
\begin{equation*}
\left\Vert d_{A}-d_{A_{0}}\right\Vert =\left\Vert L_{a}\right\Vert \leq
C\left\Vert a\right\Vert <C\varepsilon.
\end{equation*}%
Here we used \emph{Lemma %
\ref{Lemma Sobolev's embedding}}, since our choice of $\left\Vert.  \right\Vert^p_k$ suits Sobolev's embedding theorem. So $\left( \ref{restriction of dA}\right) $
is also an isomorphism for $\varepsilon $ small enough, and we can
find a unique $\rho \in \ker d_{A_{0}}^{\ast }$ solving $\left( \ref%
{surjectivity}\right) $.

\noindent\emph{\textbf{Injectivity}} 

Let $\eta \in \ker d_{A}\subset \Omega^6\left( 
\mathfrak{g}\right) $; then
\begin{eqnarray*}
        p_{a}\left( \eta \right) =0
        &\Leftrightarrow &\rho =\eta \in \ker d_{A} \\
        &\Leftrightarrow &d_{A}\rho =0 \\
        &\Leftrightarrow &\rho =0
\end{eqnarray*}%
since $\rho \in \func{img}d_{A_{0}}^{\ast }$ and we have just seen that $d_{A}:\img%
d_{A_{0}}^{\ast }\tilde{\rightarrow}\; \img d_{A_{0}}$ is an
isomorphism (for suitably small $a$); so $\eta =\rho =0$.
\end{proof}
\end{proposition}

  Hence \emph{Proposition $\ref{Prop local description of ME}$} is proved, in the terms of \emph{Corollary \ref{cor zero set}} [cf. \emph{Appendix}].

\newpage
\subsection{Final comments: gluing families and transversality}

We achieved in \emph{Theorem \ref{thm gluing}} this paper's main goal of constructing a solution $A$ of the instanton equation over the compact $G_2-$manifold $M_S=W'\tsum W''$. To conclude, I will briefly outline two natural extensions of this theory.  

First, one may consider instanton families, i.e., given (pre)compact sets  
$N\ii\subset\mc{M}^+_{\mc{E}\ii}$ of regular points on the moduli spaces over each end, define - for large neck length $S$ - an operation
\begin{displaymath}
        \tau_S:N'\times N''\rightarrow\mc{M}^+_{\mc{E}_S}.
\end{displaymath} 
This should be a diffeomorphism over its image, consisting itself of regular points. Moreover, given an adequate notion of `distance' between a connection $A^{\tsum}$ on 
$\mc{E}_S$ and $A=A'\tsum A''$, any `nearby' instanton is also obtained from such a sum. Namely, for a given integer $q$, one defines
\begin{equation}
        \ell^q_S\left(A^{\tsum} ;A',A''\right)\doteq
        \inf_{g\in\mc{G}'}\Vert g.A^{\tsum}-A' \Vert_{L^q(W'_S)}+
        \inf_{g\in\mc{G}''}\Vert g.A^{\tsum}-A'' \Vert_{L^q(W''_S)}
\end{equation}
where $W\ii_S$ are the truncations of $W\ii$ at `length' $S$. Then one should expect the following to hold:
\begin{conjecture}
Let $N\ii\subset\mc{M}^+_{\mc{E}\ii}$ be compact sets of regular points in the moduli spaces of $G_2-$instantons over $W\ii\times S^1$; for sufficiently small $\delta>0$ and large $S>0$, there are neighbourhoods $V\ii\supset N\ii$ and a smooth map
\begin{displaymath}
        \tau_S:V'\times V''\rightarrow\mc{M}^+_{\mc{E}_S}
\end{displaymath} 
such that, for appropriate choices of $q>0$ (independent of $\delta$ and $S$),\begin{enumerate}
        \item 
        $\tau_S$ is a diffeomorphism to its image, which consists of regular points;

        \item
        $\ell^q_S\left(\tau_S(A',A'');A',A''\right)\leq\delta $, 
        \quad $\forall A\ii\in V\ii$;
                 
        \item
        any instanton $A\in\mc{M}^+_{\mc{E}_S}$ such that $\ell^q_S\left(\tau_S(A',A'');A',A''\right)\leq\delta $ for some $A\ii\in V\ii$ lies in the image of $\tau_S(V'\times V'')$.
\end{enumerate}
\end{conjecture}

Finally, under a generic Morse-Bott assumption, one may envisage relaxing the non-degeneracy requirement and deal with reducible connections via the `gluing parameter' over the divisor at  infinity. If the images of  the obvious restriction maps as $r\ii:\mc{M}^+_{\mc{E}\ii}\rightarrow \mc{M}_D$ meet transversally, one would expect the moduli space $\mc{M}^+_{\mc{E}_S}$ over the compact base $M_S$ to be smoothly modelled on the product
\begin{displaymath}
        \mc{M}^+_{\mc{E}'}\times_{\mc{M}_D}\mc{M}^+_{\mc{E}''}
        =\left\{ (A',A'')\mid r'(A')=r''(A'') \right\}
\end{displaymath} 
for large values of the neck length $S$. 

All of this follows strictly, of course, the general `programme' of \cite{floer}. There are essentially two reasons to expect the analogy to carry through to our $G_2$ setting: the bounded geometry of Kovalev's manofolds, which implies the  uniform (Sobolev) bounds for the relevant operator norms [cf. \emph{Lemma \ref{Lemma noncompact Sobolev}} and \emph{Proposition \ref{prop bounded multiplication}}]; and the fact that our original solutions decay exponentially in all derivatives to a HYM connection over the divisor at infinity [cf. (\ref{eq HYM solution})], so they are in $L^p_k(W)$ for any choice of $p,k$.  

%\newpage
\appendix

\section{Chern-Simons formalism under holonomy $G_2$}
\label{Subsect Chern-Simmons}

In (3+1)-dimensional gauge theory \cite[\S2.5]%
{floer}, the Chern-Simons functional is defined on $\mathcal{B=A}/\mathcal{G}$, with integer periods, its critical points being precisely the flat
connections.
A similar theory can be formulated in higher dimensions  given a suitable closed $\left(
n-3\right) -$form \cite%
{Donaldson-Thomas}\cite{Thomas}. Here, for suitable connections over a $G_2-$manifold $(M,\varphi)$, we use the Hodge-dual $%
\ast \varphi$.

Recall that the set of connections $\mathcal{A}$ is an affine space modelled on $\Omega ^{1}\left( 
\mathfrak{g}_{P}\right) $ so, fixing a reference $A_{0}\in 
\mathcal{A}$, we have $\mc{A}=A_{0}+\Omega ^{1}\left( \mathfrak{g}_{P}\right)$ and we can define 
\begin{equation*}       \index{Chern-Simons!functional}
        \vartheta \left( A\right) =\tfrac{1}{2}\int_{M}\func{tr}\left( d_{A_{0}}a\wedge         a +\frac{2}{3}a\wedge a\wedge a\right) \wedge \ast \varphi ,
\end{equation*}
fixing $\vartheta \left( A_{0}\right) =0$. Note in passing that, since only the condition $d\ast\varphi=0$ is required, the discussion extends to cases in which the $G_2-$structure $\varphi$ is not necessarily torsion-free. Moreover, the theory remains essentially unaltered if the compact base manifold is replaced by a manifold-with-boundary, say, under a setting of connections with suitable decay towards the boundary.  

The above function is
obtained by integration of the analogous $1-$form%
\begin{equation}        \label{ro ^ phi}
                        \index{Chern-Simons!1-form@$1-$form}%
        \rho \left( a\right) _{A}=\int_{M}%
        \func{tr}\left( F_{A}\wedge a\right) \wedge \ast \varphi .
\end{equation}%
We find $\vartheta $ explicitly by integrating $\rho $ over paths $A\left( t\right) =A_{0}+ta$:%
\begin{eqnarray*}
        \vartheta \left( A\right) -\vartheta \left( A_{0}\right) 
        &=&\int_{0}^{1}\rho _{A\left( t\right) }
        \left( \dot{A}\left( t\right) \right) dt \\
        &=&\int_{0}^{1}\int_{M}\func{tr}\left( \left(F_{A_{0}}+td_{A_{0}}a
        +t^{2}a\wedge a\right) \wedge a\right) \wedge \ast \varphi \\
        &=&\tfrac{1}{2}\int_{M}\func{tr}\left( d_{A_{0}}a\wedge a
        +\frac{2}{3} a\wedge a\wedge a\right) \wedge \ast \varphi .
\end{eqnarray*}

It remains to check that $\left( \ref{ro ^ phi}\right) $ is closed, so that
this doesn't depend on the path $A\left( t\right) $. Using Stokes' theorem and $d\ast \varphi =0$,
the leading term of $\rho ,$ 
\begin{equation*}
\rho \left( a\right) _{A+b}-\rho \left( a\right) _{A}=\int_{M}\limfunc{tr}%
\left( d_{A}b\wedge a\right) \wedge \ast \varphi +O\left( \left\vert
b\right\vert ^{2}\right) ,
\end{equation*}%
is indeed symmetric:%
\begin{equation*}
\int_{M}\limfunc{tr}\left( d_{A}b\wedge a-b\wedge d_{A}a\right) \wedge \ast
\varphi =\int_{M}d\left( \limfunc{tr}\left( b\wedge a\right) \wedge \ast
\varphi \right) =0.
\end{equation*}%
Hence 
\begin{equation*}
\rho \left( a\right) _{A+b}-\rho \left( a\right) _{A}=\rho \left( b\right)
_{A+a}-\rho \left( b\right) _{A}+O\left( \left\vert b\right\vert ^{2}\right),
\end{equation*}%
and we check  that 
$\rho $ is closed comparing the reciprocal Lie
derivatives on parallel vector fields $a,b$ around a point $A$:

\begin{eqnarray*}       \label{closed 1-form}
        d\rho \left( a,b\right) _{A} 
        &=&\left( \mathcal{L}_{b}\rho \left( a\right)\right) _{A}
        -\left( \mathcal{L}_{a}\rho \left( b\right) \right) _{A} \\
        &=&\lim\limits_{h\rightarrow 0}\frac{1}{h}\left\{ \left( \rho 
        \left(a\right) _{A+hb}-\rho \left( a\right) _{A}\right) 
        -\left( \rho \left(b\right) _{A+ha}
        -\rho \left( b\right) _{A}\right) \right\} \\
        &=&\lim\limits_{h\rightarrow 0}\frac{1}{h^{2}}\;
        \underset{O\left( \left\vert h\right\vert ^{3}\right) }
        {\underbrace{\left\{ \left( \rho \left( ha\right)_{A+hb}
        -\rho \left( ha\right) _{A}\right) -\left( \rho \left( hb\right)_{A+ha}
        -\rho \left( hb\right) _{A}\right) \right\} }} \\
        &=&0.
\end{eqnarray*}%
 At least locally, then, the functional $\vartheta $
descends to the orbit space $\mathcal{B}$.

To obtain the periods of $\vartheta $ under gauge action, take $g\in \mathcal{G}$ and consider a
path $\left\{ A\left( t\right) \right\} _{t\in \left[ 0,1\right] }\subset 
\mathcal{A}$ connecting $A$ to $g.A$. The natural projection then induces a bundle%
\begin{equation*}
\begin{array}{ccc}
        \mathbf{E}_{g} & \overset{\widetilde{p_{1}}}{\longrightarrow } 
        & E \\ 
        \downarrow &  & \downarrow \\ 
        M\times \left[ 0,1\right] & \overset{p_{1}}{\longrightarrow } 
        & M%
\end{array}%
\end{equation*}%
and, using $g$ to identify the fibres $\left( \mathbf{E}_{g}\right) _{0}%
\overset{g}{\simeq }\left( \mathbf{E}_{g}\right) _{1}$, we think of $%
\mathbf{E}_{g}$ as a bundle over $M\times S^{1}$. Moreover, in some trivialisation, the path $A\left( t\right) =A_{i}\left( t\right) dx^{i}$
gives a connection $\mathbf{A}=\mathbf{A}_{0}dt+\mathbf{A}_{i}dx^{i}$ on $\mathbf{E}_{g}$:  
\begin{eqnarray*}
        \left( \mathbf{A}_{0}\right) _{\left( t,p\right) } &=&0 \\
        \left( \mathbf{A}_{i}\right) _{\left( t,p\right) } 
        &=&A_{i}\left( t\right)_{p}.
\end{eqnarray*}%
The corresponding curvature $2-$form is $F_{\mathbf{A}}=\left( F_{%
\mathbf{A}}\right) _{0i}dt\wedge dx^{i}+\left( F_{\mathbf{A}}\right)
_{jk}dx^{j}\wedge dx^{k}$:%
\begin{eqnarray*}
        \left( F_{\mathbf{A}}\right) _{0i} &=&\dot{A}_{i}\left( t\right)\\
        \left( F_{\mathbf{A}}\right) _{jk} &=&\left( F_{A}\right) _{jk}.
\end{eqnarray*}
The periods of $\vartheta$ are then of the form
\begin{eqnarray*}
        \vartheta \left( g.A\right) -\vartheta \left( A\right) &=&
        \int_{0}^{1}\rho_{A\left( t\right) }\left( \dot{A}\left( t\right)         \right) dt \\
        &=&\int_{M\times \left[ 0,1\right] }\limfunc{tr}(F_{A\left( t\right)}
        \wedge\dot{A}_{i}\left( t\right) dx^{i})\wedge dt\wedge \ast \varphi\\
        &=&\int_{M\times S^{1}}\limfunc{tr}F_{\mathbf{A}}\wedge F_{\mathbf{A}}
        \wedge\ast \varphi \\
        &=&\left\langle c_{2}\left( \mathbf{E}_{g}\right) \smallsmile 
        \left[ \ast\varphi \right] ,M\times S^{1}\right\rangle .
\end{eqnarray*}%
The K\"{u}nneth formula for the cohomology of $M\times S^{1}$ gives%
\begin{equation*}
        H^{4}\left( M\times S^{1}\right) 
        =H^{4}\left( M\right) \oplus H^{3}\left(M\right) 
        \otimes \underset{\mathbb{Z}}{\underbrace{H^{1}\left( S^{1}\right)}}
\end{equation*}%
and obviously $H^{4}\left( M\right) \smallsmile \left[ \ast \varphi \right]
=0$ so, denoting $c_{2}^{\prime }\left( \mathbf{E}_{g}\right) $ the
component lying in $H^{3}\left( M\right) $ and $S_{g}=\left[ c_{2}^{\prime
}\left( \mathbf{E}_{g}\right) \right] ^{PD}$ its Poincar\'{e} dual, we are
left with%
\begin{equation*}
\vartheta \left( g.A\right) -\vartheta \left( A\right) =\left\langle \left[
\ast \varphi \right] ,S_{g}\right\rangle .
\end{equation*}
Consequently, the periods of $\vartheta $ lie in the set%
\index{Chern-Simons!periods}%
\begin{equation*}
\left\{ \left.\int_{S_{g}}\ast \varphi \;\right\vert S_{g}\in H_{4}\left( M,\mathbb{R}%
\right) \right\} .
\end{equation*}%
That may seem odd because in general this set is \emph{dense} (there is no reason to expect $\ast\varphi$ to be an integral class). Nonetheless, as
long as our interest remains in the study of the moduli space $\mc{M}^+=Z(\rho)$
of $G_{2}$-instantons, as the \emph{critical set} of $\vartheta$,
there is nothing to worry, for the gradient $\rho =d\vartheta $ is unambiguously defined
on $\mathcal{B}$.

\section{Theory of operators}
Here is a preliminary to the proofs of \emph{Proposition \ref{prop elliptic complex}}
and \emph{Proposition \ref{prop projection between kernels}}:

\begin{lemma}
For a sequence $A\overset{S}{\rightarrow }B\overset{T}{\rightarrow }C$ of linear operators (of dense domain) between Hilbert spaces, one has:%
\begin{eqnarray*}
        \ker T=\img S 
        &\Leftrightarrow& 
        \ker S^{\ast }=\img T^{\ast }.
\end{eqnarray*}

\begin{proof}
I claim $B=\img S\oplus \ker S^{\ast }$:
\begin{eqnarray*}
b\in \left( \func{img}S\right) ^{\perp }\subset B &\Leftrightarrow& \left\langle b,Sa\right\rangle =0,\quad \forall a\in A\\
&\Leftrightarrow& \left\langle S^{\ast }b,a\right\rangle =0,\quad \forall
a\in A\\
&\Leftrightarrow& b\in \ker S^{\ast }.
\end{eqnarray*}
Since $\overline{\func{Dom}(S)}=A$, $S^{\ast }$ is closed \cite[II.16]{Brezis}, so $\ker S^{\ast }\subset B$ is closed and $B=\left( \ker S^{\ast }\right) ^{\perp }\oplus
\ker S^{\ast }$.
\emph{Mutatis mutandis, }$B=\ker T\oplus \limfunc{img}T^{\ast }$, which yields the claim by uniqueness of the orthogonal complement.
\end{proof}

\index{elliptic!complex}
\begin{corollary}
\label{cor complex is elliptic iff dual is elliptic}Let $F\overset{L_{1}}{%
\rightarrow }G\overset{L_{2}}{\rightarrow }H$ be a complex of differential
operators between vector bundles with fibrewise
inner products; if the associated symbols satisfy $\sigma
\left( L_{i}^{\ast }\right) =\left( \sigma \left( L_{i}\right) \right)
^{\ast }$, then%
\begin{equation*}
F\overset{L_{1}}{\rightarrow }G\overset{L_{2}}{\rightarrow }H\text{ is
elliptic }\Leftrightarrow H\overset{L_{2}^{\ast }}{\rightarrow }G\overset{%
L_{1}^{\ast }}{\rightarrow }F\text{ is elliptic.}
\end{equation*}
\end{corollary}
\end{lemma}

Ellip´tic complexes are closely related to Fredholm theory, which provides a local model for sets of moduli given as zeroes of sections, as discussed in \emph{Subsection \ref{subsec prelim moduli theory}} and the whole of \emph{Section \ref{sect local model moduli space}}.
I adopt the notation for the Fredholm decomposition of a map \cite[4.2.5]{4-manifolds}: 
\begin{proposition}
A Fredholm map $\Xi $ from a neighbourhood of $0$ is locally
right-equivalent to a map of the form
\begin{displaymath}
\begin{array}{rcl}
        \widetilde{\Xi }:\quad U\times F & \rightarrow & V\times G \\
        \widetilde{\Xi }\left( \xi ,\eta \right) & = & 
        \left( L\left( \xi \right),\sigma\left(\xi,\eta \right)\right)\\
\end{array}
\end{displaymath}%
where $L=\left( D\Xi \right) _{0}:$ $U$ $ \tilde{\rightarrow }$ $V$ is a
linear isomorphism, $F=\ker L$ and $G=\limfunc{coker}L$ are
finite-dimensional and $\left( D\sigma \right) _{0}=0$.
\end{proposition}

\begin{corollary}       \label{cor zero set}
A neighbourhood of $0$ in $Z\left( \Xi \right) $ is
diffeomorphic to $Z\left( \nu \right) $, where 
\begin{eqnarray*}
\begin{array}{rcl}
        \nu \; : \; F &\rightarrow& G \\
        \nu \left( \eta \right) &\doteq& \sigma \left( 0,\eta \right).
\end{array}
\end{eqnarray*}
\end{corollary}

Finally, a result borrowed from \cite{Joel} is essentially the generalisation to
Banach spaces of the fact that the determinant of a linear map is continuous:

\begin{lemma}
\label{Lemma Fine}Let $D:B_{1}\rightarrow B_{2}$ be a bounded invertible
linear map of Banach spaces with bounded inverse $Q$. If $L:B_{1}\rightarrow
B_{2}$ is another linear map with%
\begin{equation*}
\left\Vert L-D\right\Vert \leq \left( 2\left\Vert Q\right\Vert \right) ^{-1},
\end{equation*}%
then $L$ is also invertible with bounded right inverse $P$ satisfying 
\begin{equation*}
\left\Vert P\right\Vert \leq 2\left\Vert Q\right\Vert .
\end{equation*}
\end{lemma}

For its application in the proof of \emph{Proposition %
\ref{prop projection between kernels}} we are also going to need the
following \emph{Lemma}, saying that the norm of the operator
`multiplication by a function' on $L^{p}$ is controlled by a suitable Sobolev norm of that
function.
\begin{lemma}           \index{norm!of multiplication operator}
                        \label{Lemma Sobolev's embedding}
On a base manifold $M=W\times S^1$, fix $f\in L_{k}^{p}\left( M\right) $ with $k\geq \dfrac{7}{p}$; then there exists a constant $c=c\left( M\right) $ such that 
\begin{equation*}
\left\Vert L_{f}\right\Vert \leq c\left\Vert f\right\Vert _{L_{k}^{p}},
\end{equation*}%
where%
\begin{equation*}
\begin{array}{r c l}
        L_{f}:L^{p}\left( M\right) &\rightarrow& L^{p}\left( M\right) \\
        g&\mapsto& f.g%
\end{array}%
.
\end{equation*}
\begin{proof}
As a direct consequence of Sobolev's embedding [\emph{Lemma \ref{Lemma noncompact Sobolev}}], one has 
\begin{displaymath}%
\left\Vert f.g\right\Vert _{L^{p}}=\left( \int_{M}\left\vert f\right\vert
^{p}.\left\vert g\right\vert^{p} d\vol \right) ^{\frac{1}{p}}\leq \left\Vert
f\right\Vert _{C^{0}}.\left\Vert g\right\Vert _{L^{p}}
\end{displaymath} 
and so%
\begin{equation*}
\left\Vert L_{f}\right\Vert =\sup\limits_{\left\Vert g\right\Vert
_{L^{p}}=1}\left\Vert f.g\right\Vert _{L^{p}}\leq \left\Vert f\right\Vert
_{C^{0}}.
\end{equation*}
\end{proof}
\end{lemma}
%+Bibliography

%-Bibliography

\end{document}